\documentclass[12pt]{article}
\usepackage{amssymb,amsmath}
\usepackage{latexsym,bm}

\usepackage{amsfonts}
\usepackage{latexsym,amsmath,amssymb,amsfonts,epsfig,psfrag,url,graphics,ifpdf,multicol}
\usepackage{eepic,color,colordvi,amscd,amsthm}
\usepackage[section]{algorithm}
\usepackage{algpseudocode}
\usepackage[numbers,sort&compress]{natbib}
\usepackage{indentfirst,graphics,epsfig}
\usepackage{graphicx}
\usepackage{graphics}
\usepackage{tikz}
\usepackage{booktabs}
\usepackage{subfigure}

\newtheorem{theo}{Theorem}

 \theoremstyle{definition}

\theoremstyle{remark}

\newcounter{casenum}[theo]

\newcounter{subcasenum}[theo]

\newcounter{claimnum}[theo]

\allowdisplaybreaks
\usepackage{indentfirst,subfig}

\begin{document}
\thispagestyle{plain}

\begin{center} {\Large On the Tur\'{a}n number of $K_m \vee C_{2k-1}$
}
\end{center}
\pagestyle{plain}
\begin{center}
{
  {\small  Jingru Yan \footnote{ E-mail address: mathyjr@163.com}}\\[3mm]
  {\small  Department of Mathematics, East China Normal University, Shanghai 200241,  China }\\

}
\end{center}

\begin{center}

\begin{minipage}{140mm}
\begin{center}
{\bf Abstract}
\end{center}
{\small   Given a graph $H$ and a positive integer $n$, the Tur\'{a}n number of $H$ for the order $n$, denoted $ex(n,H)$, is the maximum size of a simple graph of order $n$ not containing $H$ as a subgraph. Given graphs $G$ and $H$, the notation $G \vee H$ means the joint of $G$ and $H$. $\chi(G)$ denotes the chromatic number of a graph $G$. Since $\chi(K_m \vee C_{2k-1})=m+3$ and there is an edge $e\in E(K_m \vee C_{2k-1})$ such that $\chi(K_m \vee C_{2k-1}-e)= m+2$, by the Simonovits theorem, $ex(n, K_m \vee C_{2k-1}) = \lfloor \frac{(m+1)n^2}{2(m+2)}\rfloor$ for sufficiently large $n$. In this paper, we prove that $2(m+2)k-3(m+2)-1$ is large enough for $n$.

}

\end{minipage}
\end{center}

\section{Introduction}
We will deal with only finite nontrivial simple graphs. Let $G$ be a graph with vertex set $V(G)$ and edge set $E(G)$. The order $n$ and size $e(G)$ of a graph $G$ are its number of vertices and edges, respectively. The degree of $v\in V(G)$ is $d_G(v)$. $\Delta(G)$ and $\chi(G)$ denote the maximum degree and chromatic number of a graph $G$, respectively. For terminology and notations not explicitly described in this paper, the reader is referred the books \cite{BM,WE}.

Given a graph $H$ and a positive integer $n$, the \emph{Tur\'{a}n number} of $H$ for the order $n$, denoted $ex(n,H)$, is the maximum size of a simple graph of order $n$ not containing $H$ as a subgraph. A graph is \emph{$H$-free} if it contains no $H$ as a subgraph. A graph $G$ of order $n$ is an \emph{extremal graph} for $H$ if $G$ is $H$‐free with $ex(n,H)$ edges. Determining the Tur\'{a}n number of various graphs $H$ is one of the main topics in extremal graph theory \cite{B}. A classic result of Mantel \cite{M} from 1907 states that $ex(n,K_3)= \lfloor n^2/4\rfloor$ and the balanced complete bipartite graph $K_{\lfloor n/2\rfloor,\lceil n/2 \rceil}$ is the unique extremal graph. The famous Tur\'{a}n theorem \cite{T} states that $ex(n,K_{r+1})= \lfloor (r-1)n^2/(2r)\rfloor$ and the balanced complete $r$-partite graph is the unique extremal graph.

Given graphs $G$ and $H$, the notation $G \vee H$ means the \emph{joint} of $G$ and $H$. $K_n$ and $C_n$ stand for the complete graph and cycle of order $n$, respectively. A \emph{wheel} $W_n$ is a graph obtained by connecting a single vertex to all vertices of $C_{n-1}$, that is $W_n=K_1\vee C_{n-1}$. A \emph{generalized wheel} $W_{s,t}=K_s\vee C_t$. In recent years, several authors reported results on wheels \cite{BJJ,D,DJ,Y}, generalized wheels \cite{W,WC} and disjoint union of wheels \cite{XZ}. In this paper, we consider the Tur\'{a}n number of the generalized wheel $W_{m,2k-1}=K_m \vee C_{2k-1}$.

\begin{theo}\label{th1}\cite{S}
Let $H$ be a graph with $\chi(H)\geq r+1$. If there is an edge $e$ such that $\chi(H-e)=r$, then there exists an $n_0$ such that if $n>n_0$ then the balanced complete $r$-partite graph is the unique extremal graph for $H$.
\end{theo}

Since $\chi(W_{m,2k-1})=m+3$ and there is an edge $e\in E(W_{m,2k-1})$ such that $\chi(W_{m,2k-1}-e)= m+2$, by Theorem \ref{th1} and the Tur\'{a}n theorem,
$$ex(n,W_{m,2k-1})=ex(n,K_{m+3})=\lfloor \frac{(m+1)n^2}{2(m+2)}\rfloor$$ for sufficiently large $n$. In this paper, we prove that $2(m+2)k-3(m+2)-1$ is large enough for $n$.

\begin{theo}\label{th2}\cite{D}
For all $k\geq3$ and $n\geq 6k-10$, $ex(n,W_{2k}) =\lfloor n^2/3\rfloor$.
\end{theo}

Note that when $k=2$, $K_m\vee C_{2k-1}=K_{m+3}$. Thus it suffices to consider the case $k\geq 3$.

\begin{theo}\label{th3}
If $k\geq 3$ and $n\geq 2(m+2)k-3(m+2)-1$, then
$$ex(n, W_{m,2k-1}) = \lfloor \frac{(m+1)n^2}{2(m+2)}\rfloor.$$
\end{theo}

\section{Proof of Theorem 3}

\begin{proof}
First note that the case $m=1$ of Theorem \ref{th3} is exactly Theorem \ref{th2}. Thus it suffices to consider the case $m\geq 2$.

It is easy to check that the balanced complete $(m+2)$-partite graph of order $n\geq 2(m+2)k-3(m+2)-1$ does not contain $W_{m,2k-1}$ as a subgraph, where $m\geq 2$ and $k\geq 3$. Then we have
$$ex(n,W_{m,2k-1})\geq \lfloor\frac{(m+1)n^2}{2(m+2)}\rfloor$$ for $n\geq 2(m+2)k-3(m+2)-1$.

Let $G$ be a $W_{m,2k-1}$-free graph with order $n\geq 2(m+2)k-3(m+2)-1$ and size $\lfloor\frac{(m+1)n^2}{2(m+2)}\rfloor+1$. We use induction on $m$.

Next we show that the theorem holds for $m=2$. Let $G_1$ be a $W_{2,2k-1}$-free graph with order $n\geq 8k-13$ and size $\lfloor 3n^2/8\rfloor+1$. Since
$$\frac{(n-\lfloor \frac{n}{4}\rfloor-1)n}{2}<\lfloor \frac{3n^2}{8}\rfloor+1,$$
we have $\Delta (G_1)\geq n-\lfloor n/4\rfloor$.

Suppose that $\Delta (G_1)=n-1$. Let $v\in V(G_1)$ with degree $d_{G_1}(v)=n-1$. Then $e(G_1-v)=\lfloor 3n^2/8\rfloor+1-(n-1)$. Note that $n-1\geq 8k-14> 6k-10$ for $k\geq 3$ and
$$e(G_1-v)=\lfloor \frac{3n^2}{8}\rfloor+1-(n-1)> \lfloor \frac{(n-1)^2}{3}\rfloor=ex(n-1,W_{2k}).$$ By Theorem \ref{th2}, $G_1-v$ contains $W_{2k}$ and then $G_1$ contains $W_{2,2k-1}$, a contradiction.

Next suppose $\Delta (G_1)=n-p$, where $2\leq p\leq \lfloor n/4\rfloor$. Let $v\in V(G_1)$ with degree $d_{G_1}(v)=n-p$. Let $G'_1$ be a graph of order $n-p$, obtained from $G_1$ by deleting $v$ and $p-1$ vertices not adjacent to $v$. Then
$$e(G'_1)\geq e(G_1)-p(n-p)=\lfloor \frac{3n^2}{8}\rfloor+1-p(n-p).$$
A simple calculation shows that the maximum value of the expression $p(n-p)+\lfloor (n-p)^2/3\rfloor$ when $2\leq p\leq \lfloor n/4\rfloor$ is equal to $\lfloor 3n^2/8\rfloor$. Note that $n-p\geq n-\lfloor n/4\rfloor\geq 6k-10$ for $n\geq 8k-13$ and
$$e(G'_1)\geq \lfloor \frac{3n^2}{8}\rfloor+1-p(n-p)>\lfloor \frac{(n-p)^2}{3}\rfloor=ex(n-p,W_{2k}).$$
By Theorem \ref{th2}, $G'_1$ contains $W_{2k}$ and then $G_1$ contains $W_{2,2k-1}$, a contradiction. This completes the proof of the case $m=2$ of Theorem \ref{th3}.

Therefore, by the induction hypothesis, if $n\geq 2(m+1)k-3(m+1)-1$, then $$ex(n, W_{m-1,2k-1}) = \lfloor \frac{mn^2}{2(m+1)}\rfloor.$$ Now we show that the theorem holds for $m$. Recall that $G$ is a $W_{m,2k-1}$-free graph with order $n\geq 2(m+2)k-3(m+2)-1$ and size $\lfloor\frac{(m+1)n^2}{2(m+2)}\rfloor+1$. Since $$\frac{(n-\lfloor\frac{n}{m+2}\rfloor-1)n}{2}<\lfloor\frac{(m+1)n^2}{2(m+2)}\rfloor+1,$$ we have $\Delta (G)\geq n-\lfloor n/(m+2)\rfloor$.

Suppose that $\Delta (G)=n-1$. Let $v\in V(G)$ with degree $d_{G}(v)=n-1$. Then $e(G-v)=\lfloor\frac{(m+1)n^2}{2(m+2)}\rfloor+1-(n-1)$. Note that $n-1\geq 2(m+2)k-3(m+2)-1-1> 2(m+1)k-3(m+1)-1$ for $k\geq 3$ and $$e(G-v)=\lfloor \frac{(m+1)n^2}{2(m+2)}\rfloor+1-(n-1)> \lfloor\frac{m(n-1)^2}{2(m+1)}\rfloor=ex(n-1,W_{m-1,2k-1}).$$ By Theorem \ref{th2}, $G-v$ contains $W_{m-1,2k-1}$ and then $G$ contains $W_{m,2k-1}$, a contradiction.

Next suppose $\Delta (G)=n-p$, where $2\leq p\leq \lfloor n/(m+2)\rfloor$. Let $v\in V(G)$ with degree $d_{G}(v)=n-p$. Let $G'$ be a graph of order $n-p$, obtained from $G$ by deleting $v$ and $p-1$ vertices not adjacent to $v$. Then
$$e(G')\geq e(G)-p(n-p)=\lfloor \frac{(m+1)n^2}{2(m+2)}\rfloor+1-p(n-p).$$
A simple calculation shows that the maximum value of the expression $p(n-p)+\lfloor \frac{m(n-p)^2}{2(m+1)}\rfloor$ when $2\leq p\leq \lfloor n/(m+2)\rfloor$ is equal to $\lfloor\frac{(m+1)n^2}{2(m+2)}\rfloor$. Note that $n-p\geq n-\lfloor n/(m+2)\rfloor\geq 2(m+1)k-3(m+1)-1$ for $n\geq 2(m+2)k-3(m+2)-1$ and
$$e(G')\geq\lfloor \frac{(m+1)n^2}{2(m+2)}\rfloor+1-p(n-p)>\lfloor\frac{m(n-p)^2}{2(m+1)}\rfloor=ex(n-p,W_{m-1,2k-1}).$$
By Theorem \ref{th2}, $G'$ contains $W_{m-1,2k-1}$ and then $G$ contains $W_{m,2k-1}$, a contradiction. This completes the proof of Theorem \ref{th3}.
\end{proof}

{\bf Acknowledgement} This research  was supported by the NSFC grant 12271170 and Science and Technology Commission of Shanghai Municipality (STCSM) grant 22DZ2229014.

\end{document}